\documentclass[review, 12pt]{elsarticle}
\usepackage{lineno,hyperref,amsmath}
\usepackage[noend]{algpseudocode}
\usepackage{algorithmicx,algorithm}
\usepackage{amsfonts}
\usepackage{booktabs}
\usepackage{graphics}
\usepackage{subfig}
\usepackage{float}
\usepackage{color}
\modulolinenumbers[1]
\renewcommand{\eqref}[1]{\textup{{\normalfont Eq.~(\ref{#1}}\normalfont)}}
\journal{Journal} %

\begin{document}

\begin{frontmatter}

\title{A sample-based spectral method approach for solving high-dimensional stochastic partial differential equations}

\author[a,b]{Zhibao Zheng\corref{CorrespondingAuthor}}
\cortext[CorrespondingAuthor]{Corresponding author.}
\ead{zhibaozheng@hit.edu.cn}

\author[a,b]{Hongzhe Dai}

\author[a,b]{Yuyin Wang}

\author[a,b]{Wei Wang}

\address[a]{Key Lab of Structures Dynamic Behavior and Control, Harbin Institute of Technology, Ministry of Education, Harbin 150090, China}

\address[b]{School of Civil Engineering, Harbin Institute of Technology, Harbin 150090, China}

\begin{abstract}
In the past few decades, a growing interest has been devoted for the propagation of uncertainties through physical models governed by stochastic partial differential equations (SPDEs). Despite their success and applications, existing methods are mainly limited to low-dimensional stochastic problems due to the extreme computational costs. In this paper, on the basis of an universal construct of stochastic solutions, we explore an efficient strategy for solving high-dimensional linear and nonlinear SPDEs, where original SPDEs are transformed into deterministic PDEs and one-dimensional stochastic algebraic equations by use of Galerkin method. Deterministic PDEs are solved by existing pde technologies and one-dimensional stochastic algebraic equations are calculated by a sample-based numerical method. Since computational costs are almost insensitive to the stochastic dimensions, the proposed strategy beats the so-called \emph{Curse of Dimensionality} with great success. Results on numerical examples in high dimensions, including the linear elliptic PDE with stochastic coefficients, the nonlinear Burgers equation with stochastic force and the linear wave equation with stochastic initial values, demonstrate that the proposed strategy provides an efficient and unified framework for solving SPDEs, and is particularly appropriate for high-dimensional SPDEs of practical interests.
\end{abstract}


\end{frontmatter}


\section{Introduction}

Due to the significant development in computational hardware and scientific computing techniques, it is now possible to solve very high resolution models in various computational physics problems, ranging from fluid mechanics to nano-bio mechanics. In particular, partial differential equations and closely related approximations have become state-of-the-art \cite{courant2008methods}. On the other hand, however, the considerable influence of inherent uncertainties on system behavior has led the scientific community to recognize the importance of uncertainty quantification (UQ) to realistic physical systems. More than ever, the goal becomes to represent and propagate uncertainties from the available data to the desired results through stochastic partial differential equations (SPDEs) \cite{prevot2007concise, chow2007stochastic}. In many applications, one has to deal with a large number of uncertain parameters, thus the computation of the SPDE requires the solution of a very high dimensional problem.

Over the last few decades, there has been tremendous progress in solving SPDEs.  From some perspectives, these methods broadly speaking be divided into \emph{intrusive} and \emph{non-intrusive} ones. As the most powerful non-intrusive method, Monte Carlo simulation (MC) and its variations \cite{caflisch1998monte, robert2013monte} have been widely used for solving SPDEs. MC methods are very easy to implement by use of the already existing deterministic solvers and its convergence rate does not depend on the number of independent uncertain or random variables. However, high computational costs can not be avioded since a large number of deterministic realizations are necessary to achieve a god accuracy. Another popular non-intrusive technique is sparse grid collocation approaches \cite{ganapathysubramanian2007sparse, nobile2008sparse, ma2009adaptive, xiu2010numerical}. In this scheme, the tensor product construction of quadrature rules \cite{schwab2011sparse} is an explicit dependence on the random dimensionality, which leads that the computational complexity grows exponentially fast with respect to the dimensionality, known as the so called \emph{Curse of Dimensionality} \cite{beylkin2002numerical, chen2017beating, bui2008model}. Other non-intrusive techniques \cite{doostan2009least, doostan2013non, kumar2016efficient} are also developed for high-dimensional problems.

To improve convergence rates, an intrusive method is known as Galerkin-type stochastic finite element method (SFEM) or corresponding extensions \cite{Xiu2002The, ghanem2003stochastic, Matthies2005Galerkin, nouy2007generalized, Babu2010On, le2010spectral}, for formulating and discretizing SPDEs. It has been proven efficient both numerically and analytically on numerous problems in engineering and science \cite{le2010spectral, najm2009uncertainty}. In this method, the target function is projected onto a stochastic space spanned by (generalized) polynomial chaos basis and Galerkin projection scheme is then used to transform the original SPDE into a system of coupled deterministic equations whose size can be up to orders of magnitude larger than that of the corresponding deterministic problems. The solution of such augmented algebraic systems is still challenging due to the increased memory and computational resources required, especially for large-scale problems. 
Furthermore, the \emph{Curse of Dimensionality} arises as the number of stochastic dimensions and/or the number of expansion terms increase. Some attempts are proposed, such as sparse approximation \cite{blatman2010adaptive, doostan2011non}, model reduction \cite{rabitz1999general, bui2008model, ma2010adaptive} and proper generalized decompositions \cite{chinesta2010recent, nouy2010proper}, to improve this point. However, effective treatments of the \emph{Curse of Dimensionality} in stochastic sapces are still an open problem.

Here we develop highly efficient numerical strategies for the explicit and high precision solution of SPDEs with application to problems that involve high-dimensional uncertainties. Based on separated representations, an universal construct of stochastic solutions \cite{nouy2007generalized} to general SPDEs is firstly developed. By use of this solution construct, we further develop an unified numerical strategy for solving linear and nonlinear SPDEs, where original SPDEs are transformed into deterministic PDEs and one-dimensional stochastic algebraic equations by use of the (stochastic) Galerkin method. In this way, the deterministic analysis and stochastic analysis in the solving procedure can be implemented in their individual spaces and existing PDEs techniques are available. Another beauty is that the \emph{Curse of Dimensionality} can be circumvent to great extent since all uncertainties are embeded into one-dimensional stochastic algebraic equations and the computational cost is insensitive to stochastic dimensions. Thus, our method is computationally possible to solve very high-dimensional stochastic problems encountered in science and engineering.

\section{Methodology}
In this paper, we consider weak formulations of SPDEs are written as,
\begin{equation} \label{eq:In}
  R\left( {\frac{{{\partial ^2}u}}{{\partial {x^2}}},\frac{{\partial u}}{{\partial x}},u,x,\theta } \right) = 0
\end{equation}
where $u$ is the unknown stochastic solution, $\theta$ and $x$ denote stochastic and deterministic (including the time variable) spaces, respectively. In particular, when SPDEs are defined in a high-dimensional stochastic space, i.e. $\theta : = \left\{ {{\theta _1}, \cdots ,{\theta _M}} \right\}$ with a large value of $M$, challenges arise in the solution of \eqref{eq:In} due to the so-called \emph{Curse of Dimensionality} in stochastic sapces.

An universal solution construct of \eqref{eq:In} is first developed. Although it is very natural to express the stochastic solutions by means of random field expansions, available techniques are inactive since no priori knowledge about $u(x,\theta)$ can be used. In this case, we construct the stochastic solution in the form
\begin{equation} \label{eq:u00}
  u\left( {x,\theta } \right) = \sum\limits_{i = 1}^\infty  {{\lambda _i}\left( \theta  \right){d_i}\left( x \right)} 
\end{equation}
In practical, we can truncate it at the $k$-th term as,
\begin{equation} \label{eq:u0}
  {u_k}\left( {x,\theta } \right) = \sum\limits_{i = 1}^k {{\lambda _i}\left( \theta  \right){d_i}\left( x \right)}  = {u_{k - 1}}\left( {x,\theta } \right) + \Delta {u_k}\left( {x,\theta } \right)
\end{equation}
where $\left\{ {{\lambda _i}(\theta )} \right\}_{i = 1}^k$ are random variables, $\left\{ {{d_i}\left( x \right)} \right\}_{i = 1}^k$ are deterministic functions, $\Delta {u_k}\left( {x,\theta } \right) = {\lambda _k}\left( \theta  \right){d_k}\left( x \right)$ and they are all unknown. \eqref{eq:u00} are similar to some classical expansions, such as Karhunen-Lo\`{e}ve expansion and Polynomial Chaos expansion. Karhunen-Lo\`{e}ve expansion and Polynomial Chaos expansion are special cases of \eqref{eq:u00} and kinds of spectral method approaches, thus the expansion \eqref{eq:u00} can been considered as a extended spectral approach.

Note that, solution construct of \eqref{eq:u0} is independent of the form of \eqref{eq:In}, thus it's applicable for both linear and nonlinear SPDEs. On the other hand, \eqref{eq:u0} provides a separated form of deterministic and stochastic spaces, which is possible to determine $\left\{ {{\lambda _i}(\theta )} \right\}_{i = 1}^k$ and $\left\{ {{d_i}\left( x \right)} \right\}_{i = 1}^k$ in their individual space, respectively. Hence, one requires to seek deterministic functions $\left\{ {{d_i}\left( x \right)} \right\}_{i = 1}^k$ and corresponding random variables $\left\{ {{\lambda _i}(\theta )} \right\}_{i = 1}^k$ such that the approximate solution in \eqref{eq:u0} satisfies \eqref{eq:In}.

In \eqref{eq:u0}, neither $\left\{ {{d_i}\left( x \right)} \right\}_{i = 1}^k$ nor $\left\{ {{\lambda _i}(\theta )} \right\}_{i = 1}^k$ are known \emph{a priori}, we can successively determine these unknown couples $\left\{ {{\lambda _i}\left( \theta  \right), {d_i}\left( x \right)} \right\}$ one after another via iterative methods. From this point, we can substitute \eqref{eq:u0} into \eqref{eq:In} and consider $\Delta {u_k}\left( {x,\theta } \right)$ in \eqref{eq:u0} as the stochastic increment of solution $u(x,\theta)$. However, it's not facile to determine $\lambda _k\left( \theta  \right)$ and ${d_k}\left( x \right)$ at the same time. In order to avoid this difficulty, the Galerkin method and following iterative strategy \cite{nouy2007generalized} are adopted \footnote{Here writing $R\left( {\frac{{{\partial ^2}{u_k}}}{{\partial {x^2}}},\frac{{\partial {u_k}}}{{\partial x}},{u_k},x,\theta } \right)$ as $R\left( {{u_k},x,\theta } \right)$ is a abuse of notation.},
\begin{align} 
  &\int {\left[ {R\left( {{u_{k - 1}} + \lambda _k^*{d_k},x,\theta } \right)\lambda _k^*{d_k}} \right]d\rho }  = 0 \label{eq:Galerkin1} \\
  &\int {\left[ {R\left( {{u_{k - 1}} + {\lambda _k}d_k^*,x,\theta } \right){\lambda _k}d_k^*} \right]dx}  = 0 \label{eq:Galerkin2}
\end{align}
where $\rho(\theta)$ is the cumulative distribution function of $\theta$. For the given random variable $\lambda _k^*\left( \theta  \right)$, \eqref{eq:Galerkin1} makes use of stochastic Galerkin projection to generate a deterministic partial differential equation about ${d_k}\left( x \right)$, which can be solved by existing deterministic techniques, such as finite element method \cite{hughes2012finite, reddy2014introduction}, finite difference method \cite{strikwerda2004finite}, etc. Further, the random variable ${\lambda _k}\left( \theta  \right)$ can be subsequently updated via the similar Galerkin procedure in \eqref{eq:Galerkin2} for the known $d_k^*\left( x \right)$ determined through \eqref{eq:Galerkin1}. $\left\{ {{\lambda _k}\left( \theta  \right), {d_k}\left( x \right)} \right\}$ is computed by repeating \eqref{eq:Galerkin1} and \eqref{eq:Galerkin2} until a good accuracy is achieved.

\eqref{eq:Galerkin2} derives one-dimensional stochastic algebraic equation about ${\lambda _k}\left( \theta  \right)$ as the form
\begin{equation} \label{eq:Rand00}
  g\left( {{\lambda _k}\left( \theta  \right),\theta } \right) = 0
\end{equation}
which hinders problems in high stochstic dimensions due to the \emph{Curse of Dimensionality}. Here we develop a sample-based method to overcome this difficulty: for each realization of $\left\{ {{\theta ^{\left( n \right)}}} \right\}_{n = 1}^N$, ${{\lambda _k}\left( {{\theta ^{\left( n \right)}}} \right)}$ can be obtained by solving deterministic equations as,
\begin{equation} \label{eq:Rand0}
  g\left( {{\lambda _k}\left( {{\theta ^{\left( n \right)}}} \right),{\theta ^{\left( n \right)}}} \right) = 0,~n = 1, \cdots, N
\end{equation}
It's important to note that, by computing the random variable $\lambda_k(\theta)$ from a set of its realizations, the \emph{Curse of Dimensionality} can be circumvent to great extent because the computation in \eqref{eq:Rand0} is insensitive to the dimensions of $\theta$. Even for problems with very high stochastic dimensions, the total computational cost in \eqref{eq:Rand0} for computing $\{ {\lambda _k}({\theta ^{(n)}})\} _{n = 1}^N$ are negligible for linear cases and also very low for nonlinear cases since only $N$ one-dimensional nonlinear algebraic equations are solved \cite{stoer2013introduction}.

For practical purposes, a certain number of truncated items are retained in \eqref{eq:u0}. The truncation criterion is considered as a 'global' error. In this paper, it's defined as,
\begin{equation} \label{eq:ErrGlo}
{\varepsilon _{global}} = \frac{{\int {\left[ {u_k^2\left( {x,\theta } \right) - u_{k - 1}^2\left( {x,\theta } \right)} \right]dxd\rho \left( \theta  \right)} }}{{\int {u_k^2\left( {x,\theta } \right)dxd\rho \left( \theta  \right)} }}
\end{equation}
which measures the contribution of the $k$-th stochastic increment ${\lambda _k}\left( \theta  \right){d_k}$ to the stochastic solution $u\left( \theta  \right)$ and converges to the final solution when it achieves the required precision.

Further, each couple $\left\{ {{\lambda _k}\left( \theta  \right),{d_k}} \right\}$ is solved by repeating \eqref{eq:Galerkin1} and \eqref{eq:Galerkin2}. The stop criterion is considered as a 'local' error and defined as,
\begin{equation} \label{eq:ErrLoc0}
{\varepsilon _{local}} = {{\int {{{\left[ {{d_{k,j + 1}}\left( x \right) - {d_{k,j}}\left( x \right)} \right]}^2}dx} } \mathord{\left/ {\vphantom {{\int {{{\left[ {{d_{k,j + 1}}\left( x \right) - {d_{k,j}}\left( x \right)} \right]}^2}dx} } {\int {d_{k,j}^2\left( x \right)dx} }}} \right. \kern-\nulldelimiterspace} {\int {d_{k,j}^2\left( x \right)dx} }}
\end{equation}
In practical, we normalize ${{d_k}}$ and introduce $\int {d_k^2\left( x \right)dx}  = 1$, thus the above formula becomes,
\begin{equation} \label{eq:ErrLoc}
{\varepsilon _{local}} = 2 - 2\int {{d_{k,j + 1}}\left( x \right){d_{k,j}}\left( x \right)dx}
\end{equation}
which measures the difference between ${d_{k,j}}\left( x \right)$ and ${d_{k,j + 1}}\left( x \right)$ and the calculation is stopped when ${d_{k,j + 1}}\left( x \right)$ is almost the same as ${d_{k,j}}\left( x \right)$.

\begin{algorithm}[!h]
	\caption{}
	\label{Algorithm1}
	\begin{algorithmic}[1]
		\While {${\varepsilon _{global}} > {\varepsilon _1}$}
		\label{Step4}    
		\State initial $\lambda _k^{\left( 0 \right)}\left( \theta  \right)$;
		\label{Step8}    
		\Repeat
		\label{Step9}    
		\State compute $d_k^{\left( j \right)}$ by solving \eqref{eq:Galerkin1};
		\label{Step10}  
		\State compute $\lambda _k^{\left( j \right)}\left( \theta  \right)$ by \eqref{eq:Rand0};
		\label{Step12}  
		\Until ${\varepsilon _{local}} < {\varepsilon _2}$
		\label{Step14}  
		\State ${u_k}\left( \theta  \right) = \sum\limits_{i = 1}^{k - 1} {{\lambda _i}\left( \theta  \right){d_i} + {\lambda _k}\left( \theta  \right){d_k}},\; k \ge 2$;
		\label{Step15}  
		\EndWhile
		\State ${\bf{end}}$ ${\bf{while}}$
		\label{Step16}  
	\end{algorithmic}
\end{algorithm}

The resulting procedure for approximating the solution of \eqref{eq:In} is summarized in Algorithm \ref{Algorithm1}, which includes a double-loop iteration procedure. The inner loop, which is from step \ref{Step9} to \ref{Step14}, is used to determine the couple of $(\lambda_k(\theta), d_k)$, while the outer loop, which is from step \ref{Step4} to \ref{Step16}, corresponds to recursively building the set of couples and thereby the approximate solution $u_k(\theta)$. With an initial random variable $\lambda _k^{\left( 0 \right)}\left( \theta  \right)$ given in step \ref{Step8}, $d_k^{\left( j \right)}$ can be determined in step \ref{Step10}, where superscript $j$ represents the \emph{j}-th round of iteration. With the obtained $d_k^{\left( j \right)}$, random variable $\lambda _k^{\left( j \right)}\left( \theta  \right)$ is then updated in step \ref{Step12}. The outer-loop iteration then generates a set of couples such that the approximate solution in step \ref{Step15} satisfies \eqref{eq:In}. Iteration errors $\varepsilon _{global}$ and $\varepsilon _{local}$ are calculated by \eqref{eq:ErrGlo} and \eqref{eq:ErrLoc}, and Convergence errors $\varepsilon _1$ and $\varepsilon _2$ are required precisions.

\section{Numerical Examples}
\subsection{Elliptic SPDE}
As one of the most important PDEs, elliptic PDEs \cite{gilbarg2015elliptic} have a well-developed theory and provide steady-state solutions to hyperbolic and parabolic PDEs. They are well suited to describe steady states of practical problems and have numerous applications in mathematics and physics, such as geometry,  electrostatics, continuum mechanics, heat conduction, etc. In order to better describe and predict physical phenomenon of practical interests, uncertainties, including stochastic coefficients and stochastic forces, etc., are introduced into the elliptic PDEs \cite{Matthies2005Galerkin, Babu2010On}, normally arising a challenge in high dimensional case.
Existing methods are generally powerless to high-dimensional stochastic cases. To verify the effectiveness and accuracy of the proposed method, we consider a second-order linear elliptic stochastic partial differential equation with a stochastic coefficient $c\left( {x,y,\theta } \right)$ as,
\begin{equation}\label{E11}
    - \nabla \left( {c\left( {x,y,\theta } \right)\nabla u\left( {x,y,\theta } \right)} \right) + a\left( {x,y} \right)u\left( {x,y,\theta } \right) = f\left( {x,y} \right)
\end{equation}
on ${\mathcal D} = \left[ {0,1} \right] \times \left[ {0,1} \right]$ with Dirichlet boundary ${u_{\partial {\cal D}}}\left( {x,y} \right) = 0$, where coefficients are given by $a\left( {x,y} \right) = 8$, $f\left( {x,y} \right) = 150$ and
\begin{equation}\label{E12}
  c\left( {x,y,\theta } \right) = 50 + \frac{3}{{10}}\sum\limits_{j = 1}^M {{\xi _j}\left( \theta  \right){\nu _j}{c_j}\left( {x,y} \right)}
\end{equation}
where $\left\{ {{\xi _j}\left( \theta  \right)} \right\}_{j = 1}^M$ are independent uniform random variables on $\left[ { - 0.5,0.5} \right]$ and $\left\{ {{\nu _j},{c_j}\left( {x,y} \right)} \right\}$ satisfy, 
\begin{equation}\label{E121}
  \int_{\cal D} {{e^{ - \left| {{x_1} - {x_2}} \right| - \left| {{y_1} - {y_2}} \right|}}{c_j}\left( {{x_1},{y_1}} \right)d{x_1}d{y_1}}  = \nu _j^2{c_j}\left( {{x_2},{y_2}} \right)
\end{equation}
Substituting \eqref{E12} into \eqref{E11} and making use of the finite element method for the spatial discretization with 808 nodes and 1539 triangle elements yield,
\begin{equation}\label{E13}
  \left( {\sum\limits_{j = 0}^M {{\xi _j}\left( \theta  \right){K_j}} } \right)u\left( \theta  \right) = F
\end{equation}
where ${\xi _0}\left( \theta  \right) \equiv 1$. \eqref{E13} is the well-known stochastic finite element equation and we introduce high-dimensional stochastic spaces $\left\{ {{\xi _i}\left( \theta  \right)} \right\}_{i = 1}^M$ with large values of $M$. In order to solve \eqref{E13}, we substitute \eqref{eq:u0} into it and compute couples $\left\{ {{\lambda _k}\left( \theta  \right),{d_k}} \right\}$. 
If random variable $\lambda _k\left( \theta  \right)$ has been determined or given an initial value, $d_k$ can be obtained by use of \eqref{eq:Galerkin1},
\begin{equation}\label{IM1}
  \left[ {\sum\limits_{j = 0}^M {E\left\{ {\lambda _k^2\left( \theta  \right){\xi _j}\left( \theta  \right)} \right\}{K_j}} } \right]{d_k} = E\left\{ {{\lambda _k}\left( \theta  \right)\left[ {F - \left( {\sum\limits_{j = 0}^M {{\xi _j}\left( \theta  \right){K_j}} } \right){u_{k - 1}}\left( \theta  \right)} \right]} \right\}
\end{equation}
where $E\left\{  \cdot  \right\}$ is the expectation operator and \eqref{IM1} can be simplified and rewritten as
\begin{equation}\label{IM1_0}
  {{\tilde K}_k}{d_k} = {{\tilde F}_k}
\end{equation}
where
\begin{equation}\label{IM1_00}
  \left\{ \begin{array}{l}
  {{\tilde K}_k} = \sum\limits_{j = 0}^M {{c_{kkj}}{K_j}} \\
  {{\tilde F}_k} = E\left\{ {{\lambda _k}\left( \theta  \right)F} \right\} - \sum\limits_{j = 0}^M {\sum\limits_{i = 1}^{k - 1} {{c_{kij}}{K_j}{d_i}} } \\
  {c_{ijk}} = E\left\{ {{\lambda _i}\left( \theta  \right){\lambda _j}\left( \theta  \right){\xi _k}\left( \theta  \right)} \right\}
  \end{array} \right.
\end{equation}
The size of \eqref{IM1_0} is the same as the original stochastic finite element equation \eqref{E13}, so no additional computational burden is introduced. 
Once $d_k$ has been determined through \eqref{IM1_0}, random variable $\lambda _k\left( \theta  \right)$ can be subsequently updated via \eqref{eq:Galerkin2} as,
\begin{equation}\label{IM2}
  {\lambda _k}\left( \theta  \right) = \frac{{d_k^T\left[ {F - \left( {\sum\limits_{j = 0}^M {{\xi _j}\left( \theta  \right){K_j}} } \right){u_{k - 1}}\left( \theta  \right)} \right]}}{{\sum\limits_{j = 0}^M {{\xi _j}\left( \theta  \right)d_k^T{K_j}{d_k}} }}
\end{equation}
introducing 
\begin{equation}\label{IM2_0}
  \left\{ \begin{array}{l}
  {a_k}\left( \theta  \right) = d_k^TF - \sum\limits_{j = 0}^M {\sum\limits_{i = 1}^{k - 1} {{e_{kji}}{\lambda _i}\left( \theta  \right){\xi _j}\left( \theta  \right)} } \\
  {b_k}\left( \theta  \right) = \sum\limits_{j = 0}^M {{e_{kjk}}{\xi _j}\left( \theta  \right)} 
  \end{array} \right.
\end{equation}
wher ${e_{ijk}} = d_i^T{K_j}{d_k}$. Both ${{a_k}\left( \theta  \right)}$ and ${{b_k}\left( \theta  \right)}$ are random variables and \eqref{IM2} can be rewritten as 
\begin{equation}\label{IM2_01}
  {\lambda _k}\left( \theta  \right) = \frac{{{a_k}\left( \theta  \right)}}{{{b_k}\left( \theta  \right)}}
\end{equation}
The sample-based method \eqref{eq:Rand0} is adopt to solve \eqref{IM2_01} as
\begin{equation}\label{IM2_00}
  {\lambda _k}\left( {{\theta ^{\left( n \right)}}} \right) = \frac{{{a_k}\left( {{\theta ^{\left( n \right)}}} \right)}}{{{b_k}\left( {{\theta ^{\left( n \right)}}} \right)}},~n = 1, \cdots ,N
\end{equation}
The computational cost is negligible since only $N$ times division operations (or one time vector division) are involved. \eqref{IM2_01} is efficient even for very high stochastic dimensions since all random variables $\left\{ {{\xi _i}\left( \theta  \right)} \right\}$ are embedded in random variables ${{a_k}\left( \theta  \right)}$ and ${{b_k}\left( \theta  \right)}$.

Here $N = 1 \times 10^5$ random samples, i.e. $\left\{ {{\xi _j}\left( {{\theta ^{\left( n \right)}}} \right)} \right\}_{n = 1}^{1 \times {{10}^5}},~j = 1, \cdots ,M$, are adopt. Convergence errors $\varepsilon _{global}$, $\varepsilon _{local}$ in \eqref{eq:ErrGlo}, \eqref{eq:ErrLoc} are set as $1 \times {10^{ - 6}}$ and $1 \times {10^{ - 3}}$, respectively. A personal laptop (dual-core, Intel i7, 2.40GHz) is used to test different stochastic dimensions $M$. Table~\ref{tab1} shows computatinoal costs of different stochastic dimensions and corresponding iterative errors. Only 4 or 5 retained terms in \eqref{eq:u0} can achieve the required precision and computational costs increase as the stochastic dimensions increase, but not dramatically, which demonstrates the efficiency of the proposed method.

\begin{table}[htbp]
	\centering
	\caption{Computational costs of stochastic dimensions 100 to 5000 and corresponding convergence errors in iterative processes.}
	\begin{tabular}{lrrrrrrrr}
		\toprule
		& \multicolumn{5}{c}{Iterative errors at the $k$-th retained item} & \\
		\cmidrule{2-6}
		M & k=1 & k=2 & k=3 & k=4 & k=5 & Time (s)\\
		\midrule
		100   &1&8.42e-5&4.13e-6&2.46e-7&            &3.49    \\
		1000 &1&8.06e-5&4.56e-5&2.63e-7&            &44.20  \\
		2000 &1&9.49e-5&7.51e-5&2.89e-7&            &113.43\\
		3000 &1&1.39e-4&7.96e-5&1.92e-6&3.11e-7&149.20\\
		4000 &1&1.86e-4&8.10e-5&3.95e-6&2.76e-7&182.82\\
		5000 &1&2.37e-4&8.28e-5&7.66e-6&2.73e-7&225.66\\
		\bottomrule
	\end{tabular}
	\label{tab1}
\end{table}

\begin{figure*}[htbp]
	\begin{minipage}[t]{1.0\linewidth}
		\centering
		\includegraphics[width = 1.0\textwidth]{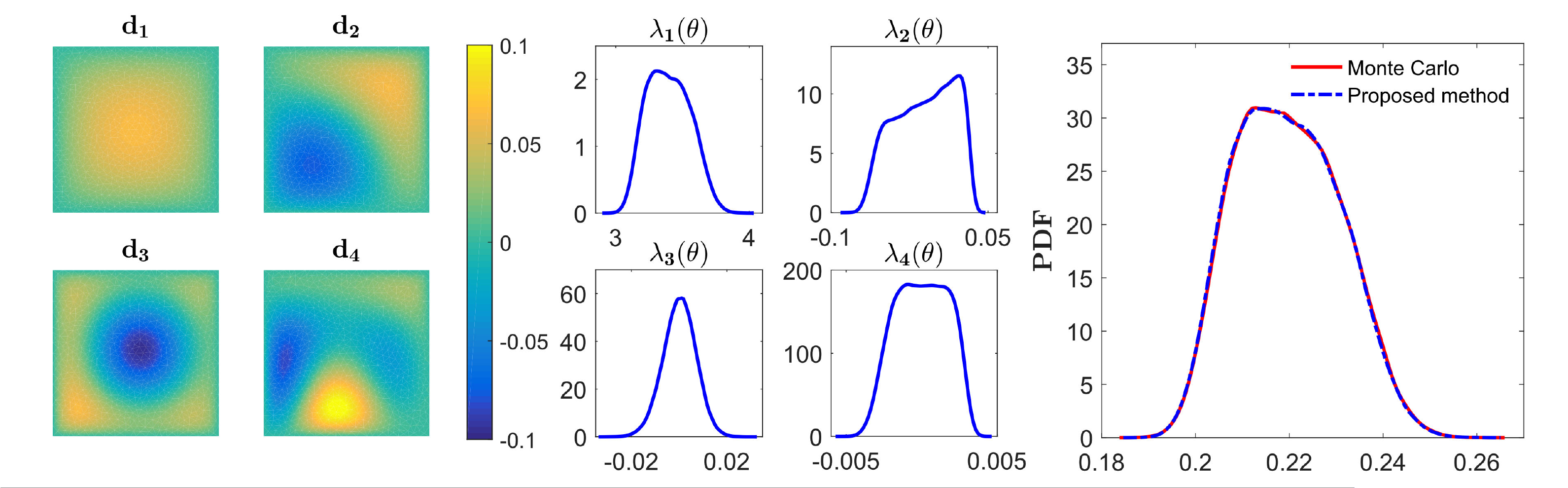}
	\end{minipage}
	\caption{Solutions $\left\{ {{d_i}\left( {x,y} \right)} \right\}_{i = 1}^4$ (left), PDFs of corresponding random variables $\left\{ {{\lambda _i}\left( \theta  \right)} \right\}_{i = 1}^4$ (mid) and comparison of PDFs at $\left( {x,y} \right) = \left( {0.5,0.5} \right)$ between $1 \times 10^6$ Monte Carlo simulations and the proposed method (right).}
	\label{Fig_E11}
\end{figure*}

\begin{figure}[!h]
	\centering
	\includegraphics[width = 0.6\linewidth]{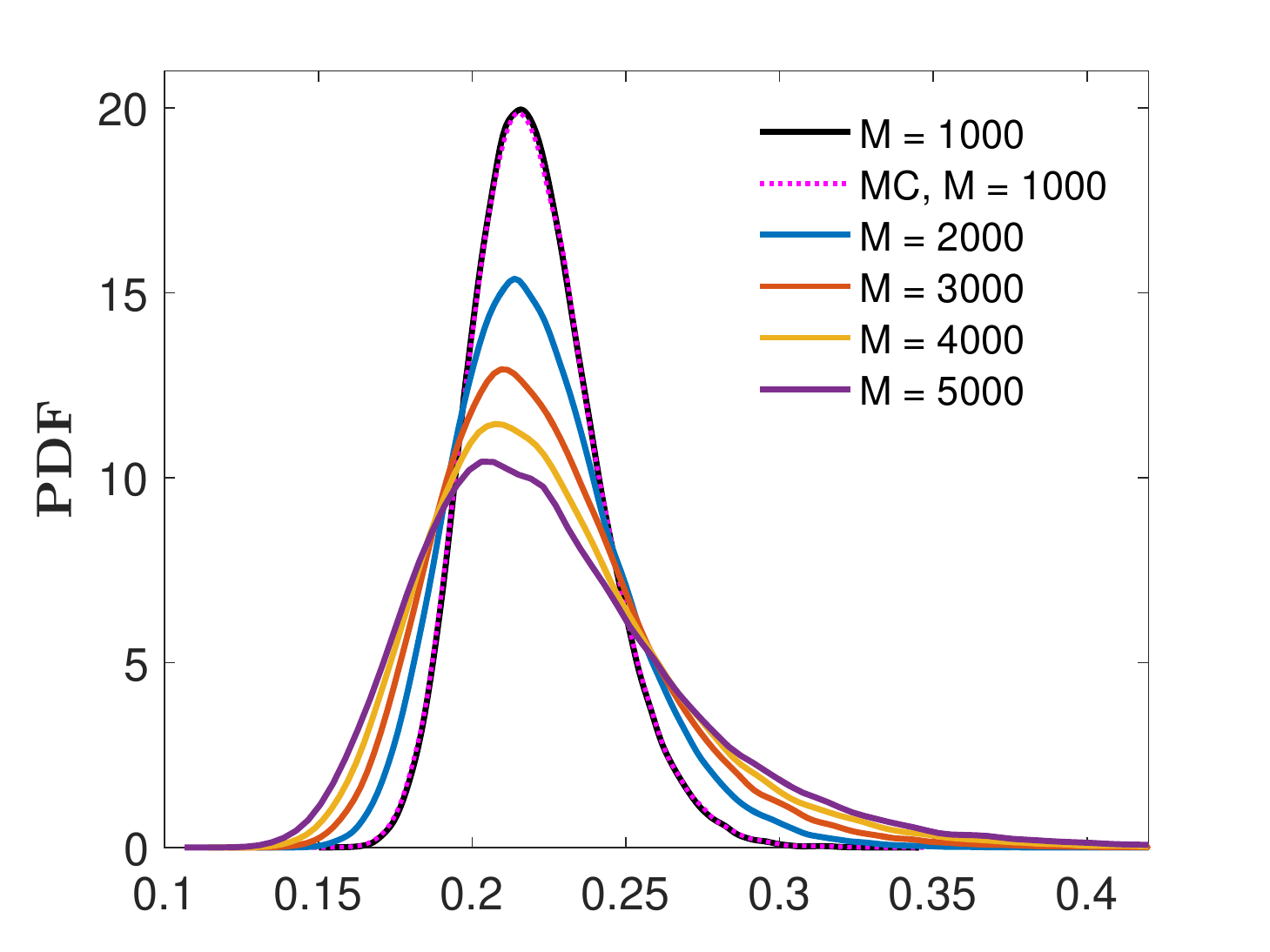}
	\caption{PDFs at $\left( {x,y} \right) = \left( {0.5,0.5} \right)$ of stochastic dimensions 1000 to 5000 and the reference solution of $M=1000$ obtained by $1 \times 10^6$ times Monte Carlo simulations.}
	\label{Fig_E12}
\end{figure}

To show some details of the proposed method, we consider the stochastic dimension $M=100$ and the reference solution is provided by $1 \times 10^6$ times Monte Carlo simulations. Figure~\ref{Fig_E11} shows solutions $\left\{ {{d_i}\left( {x,y} \right)} \right\}_{i = 1}^4$ and probability density functions (PDFs) of corresponding random variables $\left\{ {{\lambda _i}\left( \theta  \right)} \right\}_{i = 1}^4$. The comparison of PDFs between Monte Carlo simulations and the proposed method demonstrates the good accuracy of the proposed method. PDFs of high-dimensional stochastic cases are shown in Figure~\ref{Fig_E12}. Due to extreme computing costs of Monte Carlo simulations, only the PDF of $M=1000$ is compared with the reference solution, which demonstrate the good accuracy and efficiency of the proposed method for high-dimensional cases.


\subsection{Stochastic Burgers Equation}
Deterministic Burgers equation is an important PDE occurring in various areas, such as fluid mechanics, nonlinear acoustics and gas dynamics. To better model the randomness inherent of turbulence, the following stochastic Burgers equation \cite{da1994stochastic, weinan2000invariant, chorin2003averaging} with a stochastic force  is introduced,
\begin{equation}\label{E21}
  \frac{{\partial u\left( {x,t,\theta } \right)}}{{\partial t}} + \frac{1}{2}\frac{{\partial {u^2}\left( {x,t,\theta } \right)}}{{\partial x}} = \gamma \frac{{{\partial ^2}u\left( {x,t,\theta } \right)}}{{\partial {x^2}}} + f\left( {x,t,\theta } \right)
\end{equation}
on $x \times t \in \left[ { 0, 2} \right] \times \left[ {0,1} \right]$. Here we consider $\gamma \equiv 0$ and the stochastic force $f\left( {x,t,\theta } \right)$ is a Brownian motion with zero mean and covariance function $C\left( {{t_1},{t_2}} \right) = \sigma _f^2\min \left( {{t_1},{t_2}} \right)$, ${\sigma _f} = 0.2$, which can be expressed in Karhunen-Lo\`{e}ve expansion \cite{papoulis2002probability} as,
\begin{equation}\label{E22}
  f\left( {x,t,\theta } \right) = \frac{{\sqrt 2 }}{5}\sum\limits_{j = 1}^M {{\xi _j}\left( \theta  \right)\frac{{\sin \left( {j - 0.5} \right)\pi t}}{{\left( {j - 0.5} \right)\pi }}}
\end{equation}
where $\left\{ {{\xi _j}\left( \theta  \right)} \right\}_{j = 1}^M$ are independent standard gaussian random variables.

We solve \eqref{E21} by use of the central difference method and the proposed strategy, including 101 time points and 51 space nodes. Simialr to \eqref{IM1} and \eqref{IM2}, stochastic nonlinear parabolic PDE \eqref{E21} is converted into the following two equations:

a nonlinear deterministic parabolic PDE on $d_k$,
\begin{equation}\label{IM_NL1}
{h_{k1}}\frac{{\partial {d_k}}}{{\partial t}} + \frac{\partial }{{\partial x}}\left[ {{h_{k2}}d_k^2 + {h_{k3}}\left( {{u_{k - 1}}} \right){d_k}} \right] = {h_{k4}}\left( {{u_{k - 1}},f} \right)
\end{equation}
where parameters are given by,
\begin{equation*}
\left\{ \begin{array}{l}
{h_{k1}} = E\left\{ {\lambda _k^2\left( \theta  \right)} \right\}, \; {h_{k2}} = \frac{1}{2}E\left\{ {\lambda _k^3\left( \theta  \right)} \right\}\\
{h_{k3}}\left( {{u_{k - 1}}} \right) = E\left\{ {\lambda _k^2\left( \theta  \right){u_{k - 1}}\left( \theta  \right)} \right\}\\
{h_{k4}}\left( {{u_{k - 1}},f} \right) = E\left\{ {{\lambda _k}\left( \theta  \right)\left( {f - \frac{{\partial {u_{k - 1}}\left( \theta  \right)}}{{\partial t}} - \frac{1}{2}\frac{{\partial u_{k - 1}^2\left( \theta  \right)}}{{\partial x}}} \right)} \right\}
\end{array} \right.
\end{equation*}

and a one-dimensional stochastic nonlinear algebraic equation on ${\lambda _k}\left( \theta  \right)$,
\begin{equation}\label{IM_NL2}
{a_k}\lambda _k^2\left( \theta  \right) + {b_k}\left( \theta  \right){\lambda _k}\left( \theta  \right) + {c_k}\left( \theta  \right) = 0
\end{equation}
where parameters are given by,
\begin{equation*}
\left\{ \begin{array}{l}
{a_k} = \frac{1}{2}\int {{d_k}\frac{{\partial d_k^2}}{{\partial x}}dxdt} \\
{b_k}\left( \theta  \right) = \int {{d_k}\left( {\frac{{\partial {d_k}}}{{\partial t}} + \frac{{\partial \left[ {{u_{k - 1}}\left( \theta  \right){d_k}} \right]}}{{\partial x}}} \right)dxdt} \\
{c_k}\left( \theta  \right) = \int {{d_k}\left( {\frac{{\partial {u_{k - 1}}\left( \theta  \right)}}{{\partial t}} + \frac{1}{2}\frac{{\partial u_{k - 1}^2\left( \theta  \right)}}{{\partial x}} - f} \right)dxdt} 
\end{array} \right.
\end{equation*}
\eqref{IM_NL1} can be solved by the finite difference method efficiently and ${\lambda _k}\left( \theta  \right)$ in \eqref{IM_NL2} is obtained by use of the sample-based method \eqref{eq:Rand0}, which can be considered as a kind of stochastic finite difference method (SFDM). Different from \eqref{IM2}, $N$ times nonlinear algebraic equation are solved to determine $\left\{ {{\lambda _k}\left( {{\theta ^{\left( n \right)}}} \right)} \right\}_{n = 1}^N$. Computational costs increase slightly compared to the linear equation \eqref{IM2}, and it's still highly efficient for high stochastic dimensions.

\begin{figure}[hbtp]
	\centering
	\includegraphics[width = 0.6\linewidth]{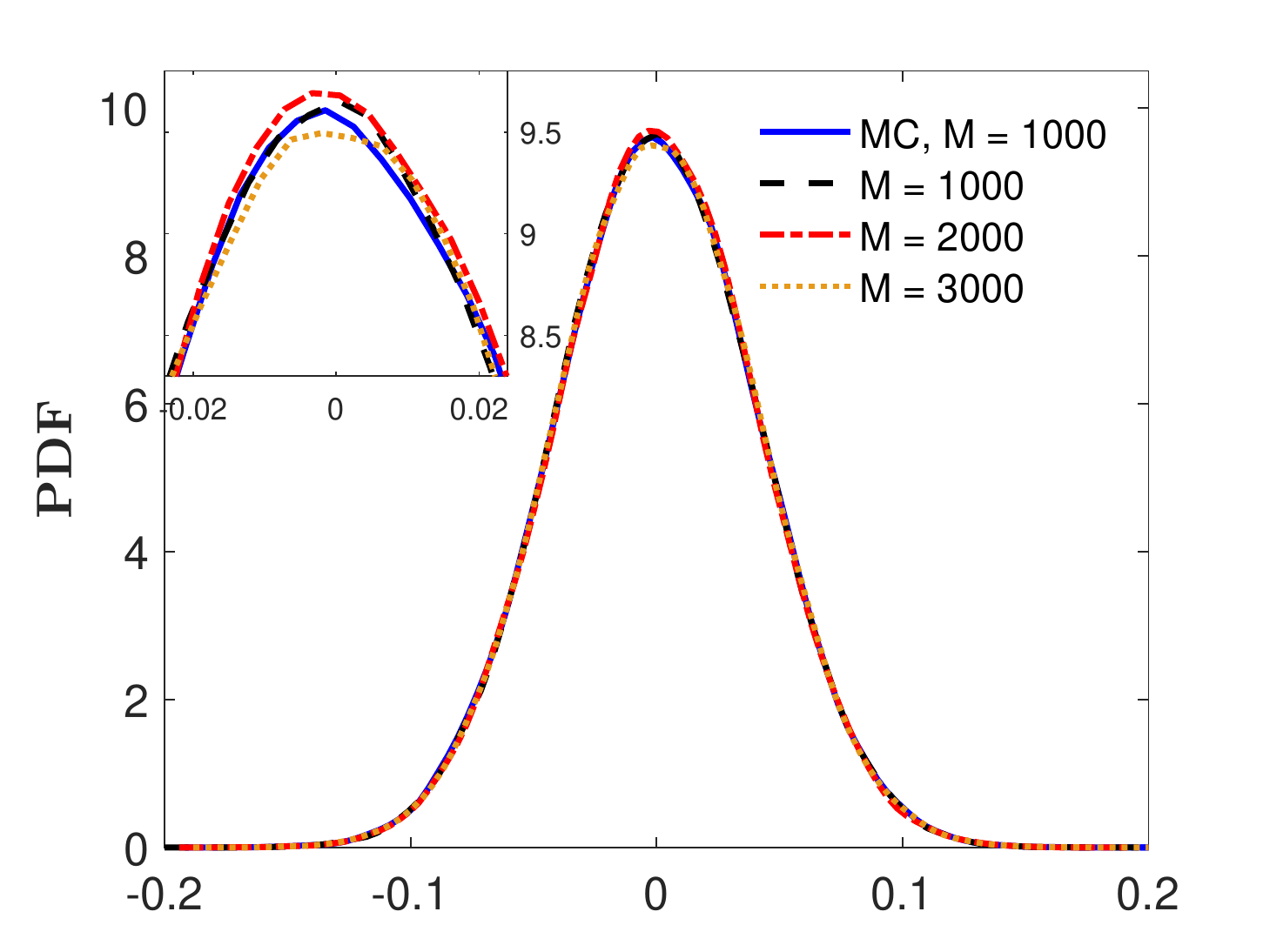}
	\caption{PDFs at $\left( {x,t} \right) = \left( {1,0.5} \right)$ of stochastic dimensions 1000, 2000, 3000 and the reference solution of $M=1000$ obtained by $1 \times 10^6$ Monte Carlo simulations.}
	\label{Fig_E21}
\end{figure}

\begin{table}[htbp]
	\centering
	\caption{Computatinoal costs of stochastic dimensions 1000, 2000, 3000 and corresponding convergence errors in iterative processes.}
	\begin{tabular}{lrrrrrr}
		\toprule
		& \multicolumn{5}{c}{Iterative errors} & \\
		\cmidrule{2-6}
		M & k=1 & k=2 & k=3 & k=4 & k=5 & Time (s)\\
		\midrule
		1000&1&1.26e-1&2.80e-2&1.11e-2&5.26e-3&159.40 \\
		2000&1&1.27e-1&2.79e-2&1.19e-2&3.98e-3&266.18 \\
		3000&1&1.26e-1&2.83e-2&1.11e-2&2.60e-3&376.88 \\
		\bottomrule
	\end{tabular}
	\label{tab2}
\end{table}
Here $N = 1 \times 10^5$ random samples and convergence errors $\varepsilon _{global} = 1 \times {10^{ - 2}}$, $\varepsilon _{local} = 1 \times {10^{ - 3}}$ are adopt. Table~\ref{tab2} shows computatinoal costs of different stochastic dimensions and corresponding iterative errors. It notes that, different from \eqref{E12}, few of retained terms $M$ in \eqref{E22} can make the stochastic soultion $u\left( {x,t,\theta } \right)$ changeless, thus iterative errors only have slight differences for large retained terms $M$. Figure~\ref{Fig_E21} shows PDFs of different stochastic dimensions and the comparison
between the reference solution and the computing solution for the stochastic dimension $M=1000$. The efficiency and accuracy are verified again.

\subsection{Stochastic Wave Equation}
As a typical representative of hyperbolic PDEs, wave equations are for the descriptions of waves occurring in many fields, such as acoustics, optics, seismology, electromagnetics and fluid dynamics \cite{courant2008methods, reed2012methods}. Here we consider a wave equation with a stochastic initial value as, 
\begin{equation}\label{E31}
  \frac{{{\partial ^2}u\left( {x,y,t,\theta } \right)}}{{\partial {t^2}}} - c\left( {x,y} \right)\Delta u\left( {x,y,t,\theta } \right) = 0
\end{equation}
on a circle with the radius 1 and $t \in \left[ {0,2} \right]$,  $c\left( {x,y} \right) = 1$. The boundary condition is ${u_{\partial {\cal D}}}\left( x,y,t \right) = 0$ and the stochastic initial value $u_{t=0}\left( {x,y,\theta } \right)$ is,
\begin{equation}\label{E32}
  {u_{t = 0}}\left( {x,y,\theta } \right) = \sqrt 2 \sum\limits_{j = 1}^M {{\xi _j}\left( \theta  \right)\sin j\pi r}
\end{equation}
and ${\left. {\frac{{\partial u\left( {x,y,t,\theta } \right)}}{{\partial t}}} \right|_{t = 0}} = 0$, where $r = \sqrt {{x^2} + {y^2}}$ denotes the polar coordinates and $\left\{ {{\xi _j}\left( \theta  \right)} \right\}_{j = 1}^M$ are independent standard gaussian random variables.

\begin{figure*}[t]
\begin{minipage}[t]{0.5\linewidth}
\centerline{\includegraphics[width=1.0\linewidth]{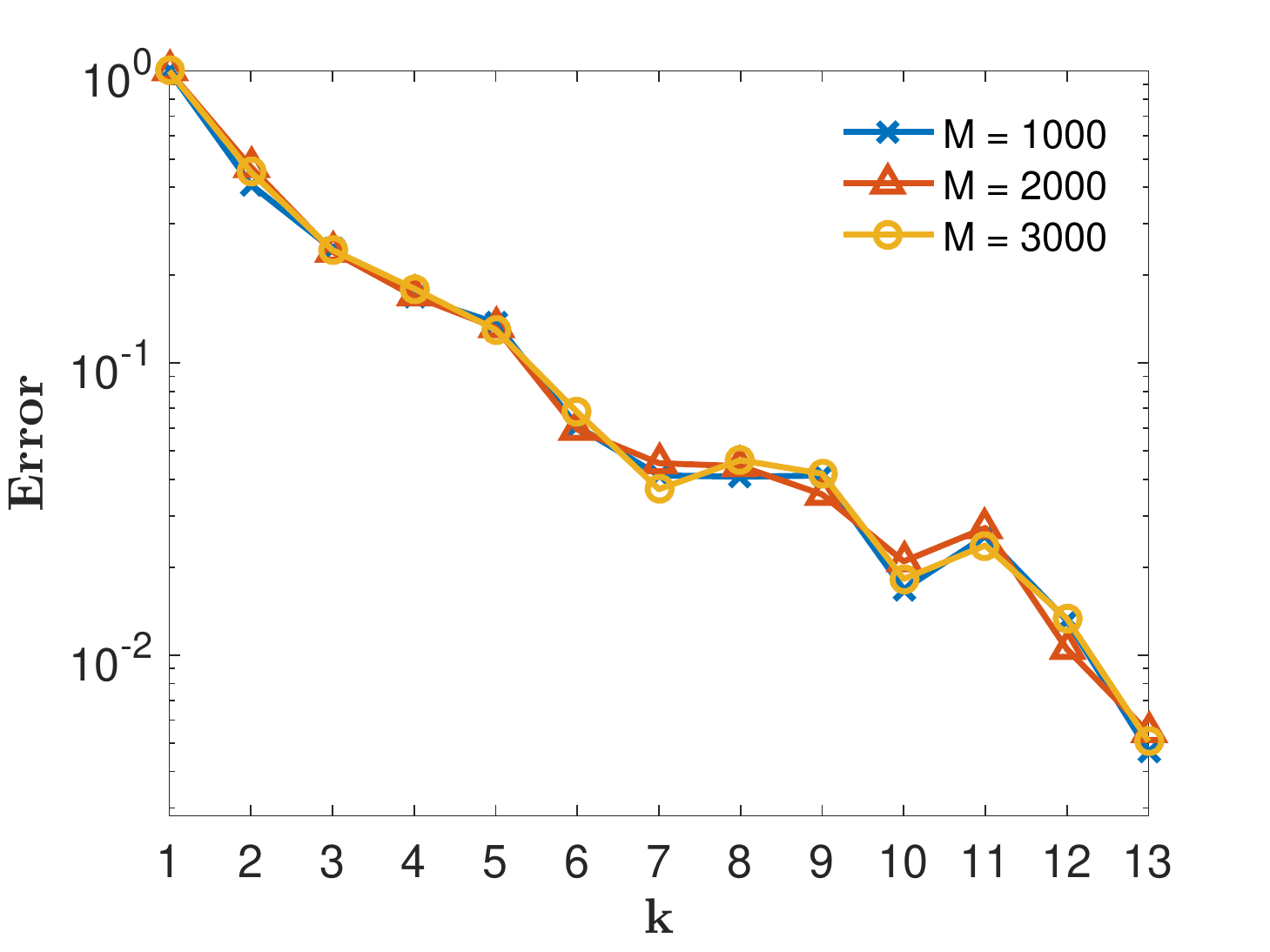}}
\end{minipage}
\begin{minipage}[t]{0.5\linewidth}
\centerline{\includegraphics[width=1.0\linewidth]{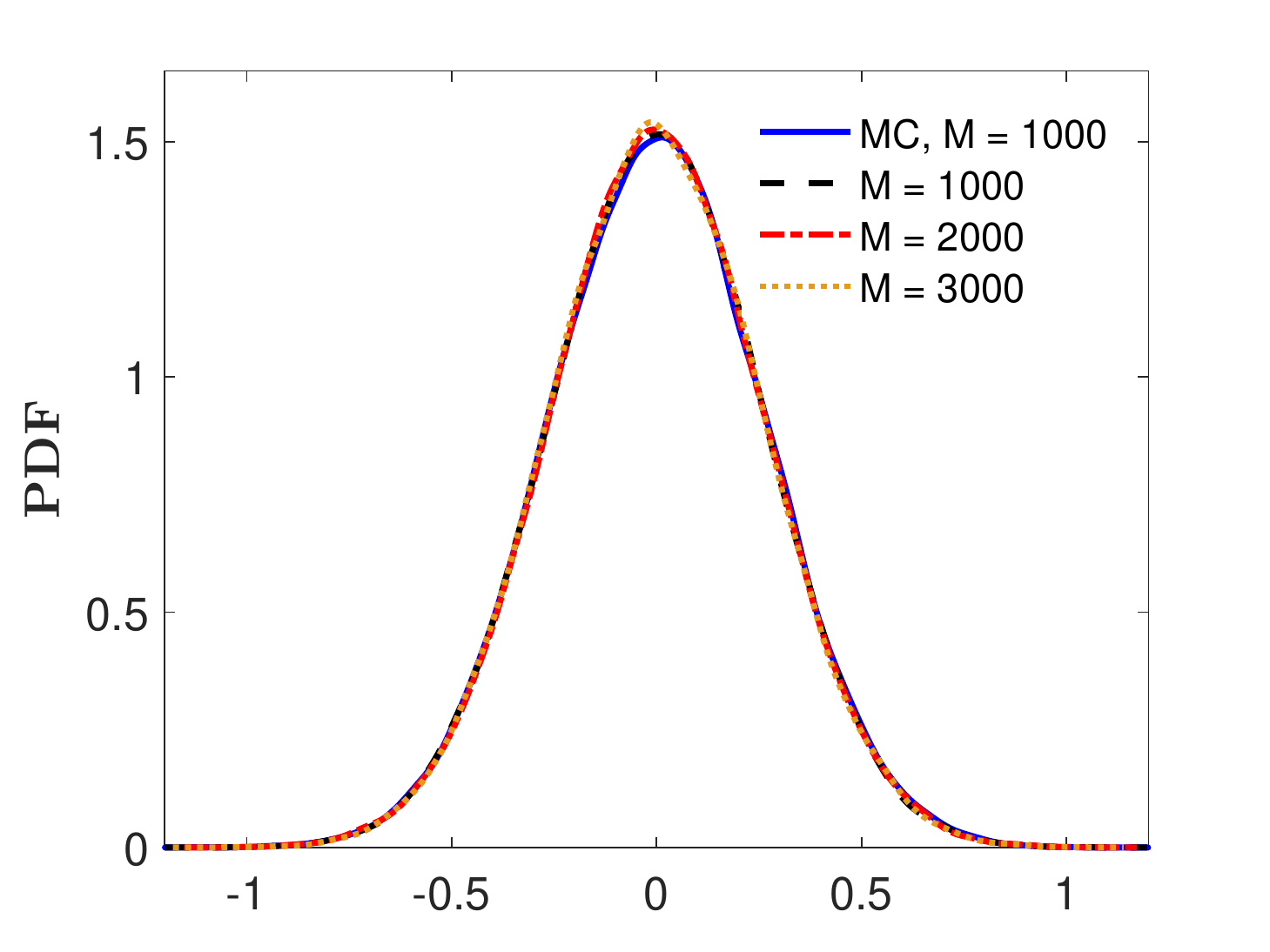}}
\end{minipage}
\caption{Convergence errors in iterative processes (left), PDFs at $\left( {x,y,t} \right) = \left( {0,0,1} \right)$ (right) of stochastic dimensions 1000, 2000 and 3000 and the reference solution of $M=1000$ obtained by $1 \times 10^6$ Monte Carlo simulations.}
\label{Fig_E31}
\end{figure*}

To solve \eqref{E31}, discretizations are achieved by use of the finite element method in space domain and the central difference method in time domain, yielding 549 nodes, 1032 triangle elements and 201 time points. $N = 1 \times 10^5$ random samples and convergence errors $\varepsilon _{global} = 1 \times {10^{ - 2}}$, $\varepsilon _{local} = 1 \times {10^{ - 3}}$ are adopt. Computational costs of different stochastic dimensions 1000, 2000 and 3000 are 162.21s, 302.53s and 424.19s, respectively. Corresponding iterative errors, PDFs and the reference solution are shown in Figure~\ref{Fig_E31}, which again demonstrate strong applicabilities of the proposed method.

\section{Conclusions}
This paper develops an efficient and unified strategy for high precision solutions of high-dimensional SPDEs, where deterministic and stochastic analysis can be implemented in individual spaces and existing analysis techniques can be readily incorporated into solving procedures. One of the most challenging issue in high-dimensional SPDEs, known as \emph{Curse of Dimensionality}, can be circumvent with great success as computational costs of the proposed method are almost insensitive to the stochastic dimensions of SPDEs. In this sense, our algorithm is general-purpose and has great potential in the uncertainty quantification in science and engineering. In the follow-up research, we hopefully further improve the theoretical analysis of proposed method \cite{bickel2018projection} and apply the method to a wider range of problems, such as complex Bayesian inference \cite{mark2018bayesian}, weather prediction \cite{alley2019advances}, etc.

\section*{Acknowledgments}
This research was supported by the National Natural Science Foundation of China (Project 11972009). This support is gratefully acknowledged.

\bibliography{References}

\begin{thebibliography}{42}
\expandafter\ifx\csname natexlab\endcsname\relax\def\natexlab#1{#1}\fi
\providecommand{\url}[1]{\texttt{#1}}
\providecommand{\href}[2]{#2}
\providecommand{\path}[1]{#1}
\providecommand{\DOIprefix}{doi:}
\providecommand{\ArXivprefix}{arXiv:}
\providecommand{\URLprefix}{URL: }
\providecommand{\Pubmedprefix}{pmid:}
\providecommand{\doi}[1]{\href{http://dx.doi.org/#1}{\path{#1}}}
\providecommand{\Pubmed}[1]{\href{pmid:#1}{\path{#1}}}
\providecommand{\bibinfo}[2]{#2}
\ifx\xfnm\relax \def\xfnm[#1]{\unskip,\space#1}\fi
\bibitem[{Courant and Hilbert(2008)}]{courant2008methods}
\bibinfo{author}{R.~Courant}, \bibinfo{author}{D.~Hilbert},
  \bibinfo{title}{Methods of Mathematical Physics: Partial Differential
  Equations}, \bibinfo{publisher}{John Wiley \& Sons}, \bibinfo{year}{2008}.
\bibitem[{Pr{\'e}v{\^o}t and R{\"o}ckner(2007)}]{prevot2007concise}
\bibinfo{author}{C.~Pr{\'e}v{\^o}t}, \bibinfo{author}{M.~R{\"o}ckner},
  \bibinfo{title}{A concise course on stochastic partial differential
  equations}, volume \bibinfo{volume}{1905}, \bibinfo{publisher}{Springer},
  \bibinfo{year}{2007}.
\bibitem[{Chow(2007)}]{chow2007stochastic}
\bibinfo{author}{P.-L. Chow}, \bibinfo{title}{Stochastic partial differential
  equations}, \bibinfo{publisher}{Chapman and Hall/CRC}, \bibinfo{year}{2007}.
\bibitem[{Caflisch(1998)}]{caflisch1998monte}
\bibinfo{author}{R.~E. Caflisch},
\newblock \bibinfo{title}{Monte carlo and quasi-monte carlo methods},
\newblock \bibinfo{journal}{Acta numerica} \bibinfo{volume}{7}
  (\bibinfo{year}{1998}) \bibinfo{pages}{1--49}.
\bibitem[{Robert and Casella(2013)}]{robert2013monte}
\bibinfo{author}{C.~Robert}, \bibinfo{author}{G.~Casella},
  \bibinfo{title}{Monte Carlo statistical methods},
  \bibinfo{publisher}{Springer Science \& Business Media},
  \bibinfo{year}{2013}.
\bibitem[{Ganapathysubramanian and
  Zabaras(2007)}]{ganapathysubramanian2007sparse}
\bibinfo{author}{B.~Ganapathysubramanian}, \bibinfo{author}{N.~Zabaras},
\newblock \bibinfo{title}{Sparse grid collocation schemes for stochastic
  natural convection problems},
\newblock \bibinfo{journal}{Journal of Computational Physics}
  \bibinfo{volume}{225} (\bibinfo{year}{2007}) \bibinfo{pages}{652--685}.
\bibitem[{Nobile et~al.(2008)Nobile, Tempone, and Webster}]{nobile2008sparse}
\bibinfo{author}{F.~Nobile}, \bibinfo{author}{R.~Tempone},
  \bibinfo{author}{C.~G. Webster},
\newblock \bibinfo{title}{A sparse grid stochastic collocation method for
  partial differential equations with random input data},
\newblock \bibinfo{journal}{SIAM Journal on Numerical Analysis}
  \bibinfo{volume}{46} (\bibinfo{year}{2008}) \bibinfo{pages}{2309--2345}.
\bibitem[{Ma and Zabaras(2009)}]{ma2009adaptive}
\bibinfo{author}{X.~Ma}, \bibinfo{author}{N.~Zabaras},
\newblock \bibinfo{title}{An adaptive hierarchical sparse grid collocation
  algorithm for the solution of stochastic differential equations},
\newblock \bibinfo{journal}{Journal of Computational Physics}
  \bibinfo{volume}{228} (\bibinfo{year}{2009}) \bibinfo{pages}{3084--3113}.
\bibitem[{Xiu(2010)}]{xiu2010numerical}
\bibinfo{author}{D.~Xiu}, \bibinfo{title}{Numerical methods for stochastic
  computations: a spectral method approach}, \bibinfo{publisher}{Princeton
  university press}, \bibinfo{year}{2010}.
\bibitem[{Schwab and Gittelson(2011)}]{schwab2011sparse}
\bibinfo{author}{C.~Schwab}, \bibinfo{author}{C.~J. Gittelson},
\newblock \bibinfo{title}{Sparse tensor discretizations of high-dimensional
  parametric and stochastic pdes},
\newblock \bibinfo{journal}{Acta Numerica} \bibinfo{volume}{20}
  (\bibinfo{year}{2011}) \bibinfo{pages}{291--467}.
\bibitem[{Beylkin and Mohlenkamp(2002)}]{beylkin2002numerical}
\bibinfo{author}{G.~Beylkin}, \bibinfo{author}{M.~J. Mohlenkamp},
\newblock \bibinfo{title}{Numerical operator calculus in higher dimensions},
\newblock \bibinfo{journal}{Proceedings of the National Academy of Sciences}
  \bibinfo{volume}{99} (\bibinfo{year}{2002}) \bibinfo{pages}{10246--10251}.
\bibitem[{Chen and Majda(2017)}]{chen2017beating}
\bibinfo{author}{N.~Chen}, \bibinfo{author}{A.~J. Majda},
\newblock \bibinfo{title}{Beating the curse of dimension with accurate
  statistics for the fokker--planck equation in complex turbulent systems},
\newblock \bibinfo{journal}{Proceedings of the National Academy of Sciences}
  \bibinfo{volume}{114} (\bibinfo{year}{2017}) \bibinfo{pages}{12864--12869}.
\bibitem[{Bui-Thanh et~al.(2008)Bui-Thanh, Willcox, and Ghattas}]{bui2008model}
\bibinfo{author}{T.~Bui-Thanh}, \bibinfo{author}{K.~Willcox},
  \bibinfo{author}{O.~Ghattas},
\newblock \bibinfo{title}{Model reduction for large-scale systems with
  high-dimensional parametric input space},
\newblock \bibinfo{journal}{SIAM Journal on Scientific Computing}
  \bibinfo{volume}{30} (\bibinfo{year}{2008}) \bibinfo{pages}{3270--3288}.
\bibitem[{Doostan and Iaccarino(2009)}]{doostan2009least}
\bibinfo{author}{A.~Doostan}, \bibinfo{author}{G.~Iaccarino},
\newblock \bibinfo{title}{A least-squares approximation of partial differential
  equations with high-dimensional random inputs},
\newblock \bibinfo{journal}{Journal of Computational Physics}
  \bibinfo{volume}{228} (\bibinfo{year}{2009}) \bibinfo{pages}{4332--4345}.
\bibitem[{Doostan et~al.(2013)Doostan, Validi, and Iaccarino}]{doostan2013non}
\bibinfo{author}{A.~Doostan}, \bibinfo{author}{A.~Validi},
  \bibinfo{author}{G.~Iaccarino},
\newblock \bibinfo{title}{Non-intrusive low-rank separated approximation of
  high-dimensional stochastic models},
\newblock \bibinfo{journal}{Computer Methods in Applied Mechanics and
  Engineering} \bibinfo{volume}{263} (\bibinfo{year}{2013})
  \bibinfo{pages}{42--55}.
\bibitem[{Kumar et~al.(2016)Kumar, Raisee, and Lacor}]{kumar2016efficient}
\bibinfo{author}{D.~Kumar}, \bibinfo{author}{M.~Raisee},
  \bibinfo{author}{C.~Lacor},
\newblock \bibinfo{title}{An efficient non-intrusive reduced basis model for
  high dimensional stochastic problems in cfd},
\newblock \bibinfo{journal}{Computers \& Fluids} \bibinfo{volume}{138}
  (\bibinfo{year}{2016}) \bibinfo{pages}{67--82}.
\bibitem[{Xiu and Karniadakis(2002)}]{Xiu2002The}
\bibinfo{author}{D.~Xiu}, \bibinfo{author}{G.~E. Karniadakis},
\newblock \bibinfo{title}{The wiener-askey polynomial chaos for stochastic
  dierential equations},
\newblock \bibinfo{journal}{Siam Journal on Scientific Computing}
  \bibinfo{volume}{24} (\bibinfo{year}{2002}) \bibinfo{pages}{619--644}.
\bibitem[{Ghanem and Spanos(2003)}]{ghanem2003stochastic}
\bibinfo{author}{R.~G. Ghanem}, \bibinfo{author}{P.~D. Spanos},
  \bibinfo{title}{Stochastic finite elements: a spectral approach},
  \bibinfo{publisher}{Courier Corporation}, \bibinfo{year}{2003}.
\bibitem[{Matthies and Keese(2005)}]{Matthies2005Galerkin}
\bibinfo{author}{H.~G. Matthies}, \bibinfo{author}{A.~Keese},
\newblock \bibinfo{title}{Galerkin methods for linear and nonlinear elliptic
  stochastic partial differential equations},
\newblock \bibinfo{journal}{Computer Methods in Applied Mechanics Engineering}
  \bibinfo{volume}{194} (\bibinfo{year}{2005}) \bibinfo{pages}{1295--1331}.
\bibitem[{Nouy(2007)}]{nouy2007generalized}
\bibinfo{author}{A.~Nouy},
\newblock \bibinfo{title}{A generalized spectral decomposition technique to
  solve a class of linear stochastic partial differential equations},
\newblock \bibinfo{journal}{Computer Methods in Applied Mechanics and
  Engineering} \bibinfo{volume}{196} (\bibinfo{year}{2007})
  \bibinfo{pages}{4521--4537}.
\bibitem[{Babu{\v{s}}ka and Chatzipantelidis(2002)}]{Babu2010On}
\bibinfo{author}{I.~Babu{\v{s}}ka}, \bibinfo{author}{P.~Chatzipantelidis},
\newblock \bibinfo{title}{On solving elliptic stochastic partial differential
  equations},
\newblock \bibinfo{journal}{Computer Methods in Applied Mechanics Engineering}
  \bibinfo{volume}{191} (\bibinfo{year}{2002}) \bibinfo{pages}{4093--4122}.
\bibitem[{Le~Ma{\^\i}tre and Knio(2010)}]{le2010spectral}
\bibinfo{author}{O.~Le~Ma{\^\i}tre}, \bibinfo{author}{O.~M. Knio},
  \bibinfo{title}{Spectral methods for uncertainty quantification: with
  applications to computational fluid dynamics}, \bibinfo{publisher}{Springer
  Science \& Business Media}, \bibinfo{year}{2010}.
\bibitem[{Najm(2009)}]{najm2009uncertainty}
\bibinfo{author}{H.~N. Najm},
\newblock \bibinfo{title}{Uncertainty quantification and polynomial chaos
  techniques in computational fluid dynamics},
\newblock \bibinfo{journal}{Annual Review of Fluid Mechanics}
  \bibinfo{volume}{41} (\bibinfo{year}{2009}) \bibinfo{pages}{35--52}.
\bibitem[{Blatman and Sudret(2010)}]{blatman2010adaptive}
\bibinfo{author}{G.~Blatman}, \bibinfo{author}{B.~Sudret},
\newblock \bibinfo{title}{An adaptive algorithm to build up sparse polynomial
  chaos expansions for stochastic finite element analysis},
\newblock \bibinfo{journal}{Probabilistic Engineering Mechanics}
  \bibinfo{volume}{25} (\bibinfo{year}{2010}) \bibinfo{pages}{183--197}.
\bibitem[{Doostan and Owhadi(2011)}]{doostan2011non}
\bibinfo{author}{A.~Doostan}, \bibinfo{author}{H.~Owhadi},
\newblock \bibinfo{title}{A non-adapted sparse approximation of pdes with
  stochastic inputs},
\newblock \bibinfo{journal}{Journal of Computational Physics}
  \bibinfo{volume}{230} (\bibinfo{year}{2011}) \bibinfo{pages}{3015--3034}.
\bibitem[{Rabitz and Ali{\c{s}}(1999)}]{rabitz1999general}
\bibinfo{author}{H.~Rabitz}, \bibinfo{author}{{\"O}.~F. Ali{\c{s}}},
\newblock \bibinfo{title}{General foundations of high-dimensional model
  representations},
\newblock \bibinfo{journal}{Journal of Mathematical Chemistry}
  \bibinfo{volume}{25} (\bibinfo{year}{1999}) \bibinfo{pages}{197--233}.
\bibitem[{Ma and Zabaras(2010)}]{ma2010adaptive}
\bibinfo{author}{X.~Ma}, \bibinfo{author}{N.~Zabaras},
\newblock \bibinfo{title}{An adaptive high-dimensional stochastic model
  representation technique for the solution of stochastic partial differential
  equations},
\newblock \bibinfo{journal}{Journal of Computational Physics}
  \bibinfo{volume}{229} (\bibinfo{year}{2010}) \bibinfo{pages}{3884--3915}.
\bibitem[{Chinesta et~al.(2010)Chinesta, Ammar, and Cueto}]{chinesta2010recent}
\bibinfo{author}{F.~Chinesta}, \bibinfo{author}{A.~Ammar},
  \bibinfo{author}{E.~Cueto},
\newblock \bibinfo{title}{Recent advances and new challenges in the use of the
  proper generalized decomposition for solving multidimensional models},
\newblock \bibinfo{journal}{Archives of Computational methods in Engineering}
  \bibinfo{volume}{17} (\bibinfo{year}{2010}) \bibinfo{pages}{327--350}.
\bibitem[{Nouy(2010)}]{nouy2010proper}
\bibinfo{author}{A.~Nouy},
\newblock \bibinfo{title}{Proper generalized decompositions and separated
  representations for the numerical solution of high dimensional stochastic
  problems},
\newblock \bibinfo{journal}{Archives of Computational Methods in Engineering}
  \bibinfo{volume}{17} (\bibinfo{year}{2010}) \bibinfo{pages}{403--434}.
\bibitem[{Hughes(2012)}]{hughes2012finite}
\bibinfo{author}{T.~J. Hughes}, \bibinfo{title}{The finite element method:
  linear static and dynamic finite element analysis},
  \bibinfo{publisher}{Courier Corporation}, \bibinfo{year}{2012}.
\bibitem[{Reddy(2014)}]{reddy2014introduction}
\bibinfo{author}{J.~N. Reddy}, \bibinfo{title}{An Introduction to Nonlinear
  Finite Element Analysis: with applications to heat transfer, fluid mechanics,
  and solid mechanics}, \bibinfo{publisher}{OUP Oxford}, \bibinfo{year}{2014}.
\bibitem[{Strikwerda(2004)}]{strikwerda2004finite}
\bibinfo{author}{J.~C. Strikwerda}, \bibinfo{title}{Finite difference schemes
  and partial differential equations}, volume~\bibinfo{volume}{88},
  \bibinfo{publisher}{Siam}, \bibinfo{year}{2004}.
\bibitem[{Stoer and Bulirsch(2013)}]{stoer2013introduction}
\bibinfo{author}{J.~Stoer}, \bibinfo{author}{R.~Bulirsch},
  \bibinfo{title}{Introduction to numerical analysis},
  volume~\bibinfo{volume}{12}, \bibinfo{publisher}{Springer Science \& Business
  Media}, \bibinfo{year}{2013}.
\bibitem[{Gilbarg and Trudinger(2015)}]{gilbarg2015elliptic}
\bibinfo{author}{D.~Gilbarg}, \bibinfo{author}{N.~S. Trudinger},
  \bibinfo{title}{Elliptic partial differential equations of second order},
  \bibinfo{publisher}{springer}, \bibinfo{year}{2015}.
\bibitem[{Da~Prato et~al.(1994)Da~Prato, Debussche, and
  Temam}]{da1994stochastic}
\bibinfo{author}{G.~Da~Prato}, \bibinfo{author}{A.~Debussche},
  \bibinfo{author}{R.~Temam},
\newblock \bibinfo{title}{Stochastic burgers' equation},
\newblock \bibinfo{journal}{Nonlinear Differential Equations and Applications
  NoDEA} \bibinfo{volume}{1} (\bibinfo{year}{1994}) \bibinfo{pages}{389--402}.
\bibitem[{Weinan et~al.(2000)Weinan, Khanin, Mazel, and
  Sinai}]{weinan2000invariant}
\bibinfo{author}{E.~Weinan}, \bibinfo{author}{K.~Khanin},
  \bibinfo{author}{A.~Mazel}, \bibinfo{author}{Y.~Sinai},
\newblock \bibinfo{title}{Invariant measure for burgers equation with
  stochastic forcing},
\newblock \bibinfo{journal}{Annals of Mathematics-Second Series}
  \bibinfo{volume}{151} (\bibinfo{year}{2000}) \bibinfo{pages}{877--960}.
\bibitem[{Chorin(2003)}]{chorin2003averaging}
\bibinfo{author}{A.~J. Chorin},
\newblock \bibinfo{title}{Averaging and renormalization for the
  korteveg--devries--burgers equation},
\newblock \bibinfo{journal}{Proceedings of the National Academy of Sciences}
  \bibinfo{volume}{100} (\bibinfo{year}{2003}) \bibinfo{pages}{9674--9679}.
\bibitem[{Papoulis and Pillai(2002)}]{papoulis2002probability}
\bibinfo{author}{A.~Papoulis}, \bibinfo{author}{S.~U. Pillai},
  \bibinfo{title}{Probability, random variables, and stochastic processes},
  \bibinfo{publisher}{Tata McGraw-Hill Education}, \bibinfo{year}{2002}.
\bibitem[{Reed(2012)}]{reed2012methods}
\bibinfo{author}{M.~Reed}, \bibinfo{title}{Methods of modern mathematical
  physics: Functional analysis}, \bibinfo{publisher}{Elsevier},
  \bibinfo{year}{2012}.
\bibitem[{Bickel et~al.(2018)Bickel, Kur, and Nadler}]{bickel2018projection}
\bibinfo{author}{P.~J. Bickel}, \bibinfo{author}{G.~Kur},
  \bibinfo{author}{B.~Nadler},
\newblock \bibinfo{title}{Projection pursuit in high dimensions},
\newblock \bibinfo{journal}{Proceedings of the National Academy of Sciences}
  \bibinfo{volume}{115} (\bibinfo{year}{2018}) \bibinfo{pages}{9151--9156}.
\bibitem[{Mark et~al.(2018)Mark, Metzner, Lautscham, Strissel, Strick, and
  Fabry}]{mark2018bayesian}
\bibinfo{author}{C.~Mark}, \bibinfo{author}{C.~Metzner},
  \bibinfo{author}{L.~Lautscham}, \bibinfo{author}{P.~L. Strissel},
  \bibinfo{author}{R.~Strick}, \bibinfo{author}{B.~Fabry},
\newblock \bibinfo{title}{Bayesian model selection for complex dynamic
  systems},
\newblock \bibinfo{journal}{Nature communications} \bibinfo{volume}{9}
  (\bibinfo{year}{2018}) \bibinfo{pages}{1803}.
\bibitem[{Alley et~al.(2019)Alley, Emanuel, and Zhang}]{alley2019advances}
\bibinfo{author}{R.~B. Alley}, \bibinfo{author}{K.~A. Emanuel},
  \bibinfo{author}{F.~Zhang},
\newblock \bibinfo{title}{Advances in weather prediction},
\newblock \bibinfo{journal}{Science} \bibinfo{volume}{363}
  (\bibinfo{year}{2019}) \bibinfo{pages}{342--344}.

\end{thebibliography}

\end{document}